\documentclass[12pt]{article}
\usepackage{graphicx}
\usepackage{epstopdf} 
\usepackage{amssymb} 
\usepackage{amsmath} 
\DeclareMathOperator{\Dim}{Dim}
\DeclareMathOperator{\Tr}{Tr}
\numberwithin{equation}{section}
\newtheorem{thm}{Theorem}[section]
\newtheorem{prop}[thm]{Proposition}

\begin{document}

 \centerline{\LARGE    Graphs and Generalized  Witt  \vspace{3mm}  identities}
 \centerline{\large G.A.T.F.da Costa
\footnote{g.costa@ufsc.br}}
\centerline{\large Departamento
de Matem\'{a}tica} 
\centerline{\large Universidade Federal de
Santa Catarina} \centerline{\large
88040-900-Florian\'{o}polis-SC-Brasil}

\begin{abstract}

This paper is about the determinantal  identities associated with the Ihara (Ih) zeta function of a non directed graph and  the Bowen-Lanford (BL) zeta function of a directed graph. They will be called the Ih and the BL identities in this paper.  We show that the Witt identity (WI) is a special case of the BL  identity and inspired by  the links the WI has with Lie algebras and combinatorics we investigate similar aspects of the
Ih and BL identities.
We show that they satisfy generalizations of the Strehl identity and Carlitz, Metropolis-Rota relations and each one of them  can be interpreted as the denominator (or generalized Witt) identity of a free Lie superalgebra. Also, they can be associated to a coloring problem. 
 New interpretations
of the Ih and BL zeta functions are presented.


\vspace{5mm}

\end{abstract}

\section{ Introduction}

First, I will define the Ihara (Ih) and the Bowen-Lanford (BL) zeta functions and their respective identities. 

The following definitions are needed. Let $G=(V,E)$ be a finite  connected {\it non directed} graph with no 1-degree vertices, $V$ is the set of
$|V|$ vertices, $E$ is the set of $|E|$ edges with elements  labeled $e_{1}$,
...,$e_{|E|}$. 
The
graph may have multiple edges and loops. In order to define a closed path in a non directed graph we shall orient the edges. The orientation is arbitrary but fixed. An oriented edge has origin and  end given by its orientation.
Let $G'$ be the graph with $2|E|$ oriented
edges built from the previous oriented graph $G$ by adding
in the opposing oriented edges $e_{|E|+1}=(e_{1})^{-1}$,
...,$e_{2 |E|}=(e_{|E|})^{-1}$, $(e_{i})^{-1}$ being the oriented edge
opposite to $e_{i}$ and with origin (end) the end (origin) of
$e_{i}$. In the case that $e_{i}$ is an oriented loop,
$e_{i+|E|}=(e_{i})^{-1}$ is just an additional oriented loop hooked
to the same vertex.

A path of length $N$ in the non oriented graph $G$ is the ordered sequence $(e_{i_{1}},...,e_{i_{N}})$
of oriented edges  in  $G'$ such that the end of $e_{i_{k}}$ is the origin of $e_{i_{k+1}}$. Sometimes, it will be useful to represent  a path  by a word in the alphabet of the symbols in the set $\{e_{1}, ..., e_{2E} \}$, a word being a concatenated product of symbols which respect the order of the symbols in the sequence.
A cycle is a non-backtracking tail-less closed path, that is, the end of $e_{i_{N}}$ coincides with the origin of $e_{i_{1}}$, subjected to the non-backtracking condition that $e_{i_{k+1}}  \neq e_{i_{k}+|E|}$. In another words,
a cycle  never goes immediately backwards over a  previous edge. Tail-less means that $e_{i_{1}}  \neq e_{i_{N}}^{-1}$.
The length of a cycle is the number of edges in its sequence. A cycle $p$ is called non periodic if it is not the repetition of some simpler cycle. It is periodic 
if $p=q^r$ for some $r>1$ and nonperiodic cycle $q$. The integer
$r$ is the period of $p$. The cycle $(e_{i_{N}}, e_{i_{1}}, ...,e_{i_{N-1}})$ is called a circular permutation of $(e_{i_{1}},...,e_{i_{N}})$
and $(e_{i_{N}}^{-1},...,e_{i_{1}}^{-1})$ is an inversion of the latter. 
A cycle and its inverse are taken as distinct paths whereas circular permutations are taken as equivalent.
In order to count cycles of a given length we use the {\it edge adjacency matrix} of $G$. This 
is the $2|E| \times 2|E|$ matrix $T$ defined on $G'$ as follows:
$T_{ij}=1$, if end vertex of edge $i$ is the start vertex of
edge $j$ and edge $j$ is not the inverse edge of $i$;
$T_{ij}=0$, otherwise. The number of equivalence classes of non periodic cycles of length $N$ in $G$ is given by
\begin{equation} \label{(2.3)}
\theta(N,T)=\frac{1}{N} \sum_{g|N} \mu(g) \Tr T^{\frac{N}{g}},
\end{equation}
where $g$ ranges over the positive divisors of $N$. The summation is over all positive divisors of $N$ and $\mu$ is the M\"obius function: $\mu(+1)=1$, $\mu(g)=0$, if $g=p_{1}^{a_{1}}...p_{q}^{a_{q}}$ with $a_{i}>1$, and $\mu(p_{1}...p_{q})=(-1)^{q}$, $p_{1}$, ..., $p_{q}$ primes.
The  Ih zeta function of the undirected graph $G$ is here formally defined as the reciprocal of the infinite product in the identity
\begin{equation} \label{(2.2)}
\prod_{N=1}^{+\infty} (1-z^{N})^{\theta(N,T)}=\det{(I-z T)}.
\end{equation}
In the present paper, relations (1.1) and (1.2) are called the Ih formula and the Ih identity, respectively.
See [23,30] and references therein for a very nice overview about the Ih function and identity and its properties.

Now, let $G=(V,E)$ be a finite  connected and directed graph with no 1-degree vertices, $V$ is the set of
$|V|$ vertices, $E$ is the set of $|E|$ edges with elements  labeled $e_{1}$,
...,$e_{|E|}$. 
The
graph may have multiple edges and loops. 
A path in a finite connected and directed graph $G$ is given by an ordered sequence $(e_{i_{1}},...,e_{i_{N}})$, $i_{k} \in \{1, ..., |E|\}$, of oriented edges  such that the end of $e_{i_{k}}$ is the origin of $e_{i_{k+1}}$. Notice that in this case and contrary to the previous one there are no inverse edges, hence,  paths are backtrack-less, tail-less and have no inverse. A path has a natural orientation which is induced by the orientations of the edges in the sequence.

We shall consider directed graphs which are strongly connected.
A directed graph  is called strongly connected  if it contains a directed path from $a$ to $b$ and a directed path from $b$ to $a$ for every pair of vertices $a, b$. In general, an oriented  graph may not be strongly connected but it may have strongly connected components which are the maximal strongly connected subgraphs. The graph become acyclic if each component is contracted to a single vertex. 

In order to count  cycles of a given length in a directed graph $G$ the  directed vertex adjacency matrix $A_{d}(G)$ can be used. See [1,26].
This is the matrix of order  $|V| \times |V|$ with entries defined as follows. Label the vertices of $G$, $1$, $2$, ..., $|V|$. Then, $(A_{d})_{ij}$ is the number of edges directed from vertex $i$ to vertex $j$. The number of cycles of length $N$ in $G$ is given by $\Tr A_{d}^{N}$.  Also, one can use the  directed edge adjacency matrix $T_{d}$ of $G$ to count cycles of a given length in a directed graph. We mention this result since it is not (well) known (or used) and to make a parallel with the undirected case where the matrix $T$ is always used. The matrix $T_{d}$
is $|E| \times |E|$ and defined as follows:
$(T_{d})_{ij}=1$, if end vertex of edge $i$ is the start vertex of
edge $j$;
$(T_{d})_{ij}=0$, otherwise. The number of cycles of length $N$ in $G$ is given by $\Tr T_{d}^{N}$, hence, $\Tr T_{d}^{N}= \Tr A_{d}^{N}$.  
The number of equivalence classes of non periodic cycles of length $N$ in the directed graph $G$ is $\theta_{d}$, given by (1.1) replacing $T$ by $A_{d}$ or $T_{d}$. The BL formula for $\theta_{d}$, then, follows.
The BL zeta function  of the directed graph $G$ is defined as the reciprocal of the infinite product in the BL identity obtained from (1.2) after replacing $\theta$ by $\theta_{d}$ and $T$ by $A_{d}$ or $T_{d}$. See [11,13,16,17,19,28,29].

In the sequel, I will recall  Witt identity (WI) and some facts about. Then, inspired by the links the WI has with Lie algebras, combinatorics and graphs we will investigate in the following sections similar aspects of the Ih and BL identities.
Let $N$ be a positive integer,  $R$  a real number. 
The polynomial of degree $N$ in $R$ with rational coefficients given in terms of M\"obius function,
\begin{equation} 
{\mathcal M}(N;R) = \frac{1}{N} \sum_{g \mid N} \mu (g) R^{\frac{N}{g}},
\end{equation}
is called  the Witt polynomial or the
 Witt dimension formula
according to the context where it appears, combinatorics or algebra. See 
[15,18,21].
Furthermore, it satisfies the WI
\begin{equation} \label{(1.2)}
\prod_{N=1}^{\infty} (1-z^{N})^{{\mathcal M}(N;R)} = 1- Rz.
\end{equation}

Relations (1.3) and (1.4) are associated with the following result [21]:

\begin{prop} 
If ${\mathcal V}$ is an $R$-dimensional vector space and $L$ is the free Lie algebra generated by  ${\mathcal V}$ then
$L= \bigoplus_{N=1}^{\infty} L_{N}$, and $L_{N}$ has dimension given by
${\mathcal M}(N;R)$. The generating function for the dimensions of the homogeneous subspaces of the enveloping algebra of $L$ is given by the
reciprocal of the Witt identity.
\end{prop}

The Witt polynomial (1.3) is also called the necklace polynomial  because ${\mathcal M}(N;R)$ gives the number of
 inequivalent non-periodic colorings of  a circular string of $N$ beads -  a necklace -  with at most $R$ colors. In reference [2] L. Carlitz proved that
 \begin{equation} 
 {\mathcal M} (N,\alpha \beta)=
\sum_{[s,t]=N} {\mathcal M} (s,\alpha) {\mathcal M} (t,\beta)
\end{equation}
as a special case of a more general result. The summation is over the set of positive integers $\{s,t \mid [s,t]=N\}$,  $[s,t]$ being the least common multiple of $s$ and $t$.
 In [15] Metropolis and Rota gave a new proof of this result and of other identities that are  satisfied by the Witt polynomials. For instance, 
\begin{equation} 
{\mathcal M} (N,{\alpha}^{l})=\sum_{[l,t]=N l} {\mathcal M} (t,\alpha),
\end{equation} 
where $l$ is a positive integer.
They gave a ring theoretical interpretation to these relations in order to obtain several results about the necklace ring and the Witt vectors.

A symmetrical form of the WI is given by the Strehl identity [26]:
\begin{equation}
\prod_{k \geq 1} \left( \frac{1}{1-\beta z^{k}} \right)^{{\mathcal M}(k,\alpha)}
=\prod_{j \geq 1} \left(
\frac{1}{1-\alpha z^{j}} \right)^{{\mathcal M}(j,\beta)}.
\end{equation} 

In [22] Sherman pointed out that ${\mathcal M}(N;R)$  is the number of equivalence classes of closed
non-periodic paths of length $N$ which traverse counterclockwise the edges of a graph
which has $R$  loops counterclockwise directed and hooked to a single vertex as $G_{1}$ in Figure 1. See [4,5] for a proof. This sugests a connection with the BL identity. Indeed,
the WI is the BL identity when the graph is as in Figure 1 (or an equivalent graph with respect to the BL zeta function. See Remark 3.5, section 3).   
In this case, $A_{d}=(R)$, $\Tr A_{d}^{N}=R^{N}$, and $\det(1-zA_{d})=1-Rz$, $T_{d}$ is the $R \times R$ matrix with all entries equal to one,  $\Tr T_{d}^{N}=R^{N}$,
 $\det(1-zT_{d})=1-Rz$, and $\theta_{d}(N,T_{d})= {\mathcal M}(N;R)$.

\begin{center}
\begin{figure}[ht]
\centering
\includegraphics[scale=0.5]{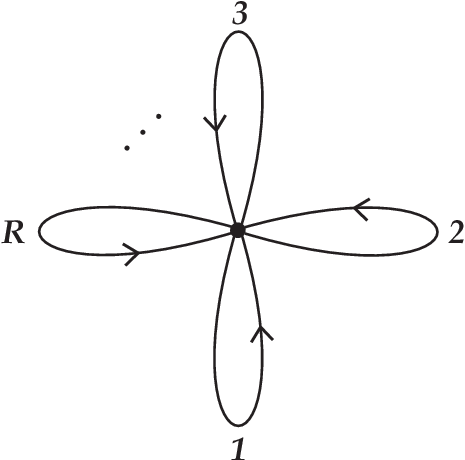}
\caption{Graph $G_{1}$}
\label{Fi:G1}
\end{figure}
\end{center}

The following definition  encapsulates the Ih and BL formulas and identities into a single formula and identity.
 
\noindent{\bf Definition.} Let  ${\mathcal T}$ stands for the  matrices $T$, $T_{d}$ or $A_{d}$ of a graph $G$ and ${\mathcal C}(N,{\mathcal T})$ stands for the respective cycle counting number formula, namely,
\begin{equation} 
{\mathcal C}(N,{\mathcal T})= \frac{1}{N}
 \sum_{g|N} \mu(g) \hspace{1mm} \Tr {\mathcal T}^{\frac{N}{g}},
\end{equation}
which satisfy the determinantal identity
\begin{equation} 
\prod_{N=1}^{+\infty} (1-z^{N})^{{\mathcal C}(N,{\mathcal T})}
=\det(I-z{\mathcal T}).
\end{equation} 
 The determinant is a polynomial of degree $D$, $D=2|E|$, if ${\mathcal T}=T$; $D=|V|$, if ${\mathcal T}=A_{d}$; and $D=|E|$, if ${\mathcal T}=T_{d}$.

The paper is organized as follows. In section 2, we show that the identity (1.9) 
and formula (1.8) satisfy generalizations of the Strehl identity and of the Carlitz, Metropolis-Rota relations. In [7-10], S. -J. Kang and collaborators generalized Proposition 1.1 to more general Lie algebras. 
In section 3, using these more general results, relations (1.8-9) are associated with  a 
free Lie superalgebra where they play the role of a dimension formula and the denominator (or generalized Witt) identity of the algebra. 
This
 establishes a connection going from graph theoretical ideas to  some  of S. -J. Kang's foundational results in [7-10] but
the results  can be understood
as a graph theoretical representation 
of some  results in [7-10]. New interpretation of  the Ih and BL zeta functions will be presented based on this connection with Lie algebras. Section 4 is devoted to examples.
In section 5, (1.8-9) are interpreted  in terms of  a restricted necklace coloring problem where another interpretation
of the Ih and BL zeta functions arises.

\section{A generalization of the Strehl identity and Carlitz, Metropolis-Rota relations} \label{sec:ids}

We start this section proving a generalization of the  Strehl identity (1.7).

\begin{thm} {\it Given two  graphs ${\mathcal G}_{1}$ and ${\mathcal G}_{2}$ with matrices ${\mathcal T}_{1}$ and ${\mathcal T}_{2}$, respectively, the following identity holds:
\begin{equation} 
\prod_{k \geq 1} \left[ \frac{1}{\det(1-z^{k}{\mathcal T}_{1})} \right]^{{\mathcal C}(k,{\mathcal T}_{2})}
=\prod_{j \geq 1} \left[
\frac{1}{\det(1-z^{j}{\mathcal T}_{2})} \right]^{{\mathcal C}(j,{\mathcal T}_{1})}.
\end{equation} }
\end{thm}

\noindent{\bf Proof.}
Call $I$ and $II$ the left and right hand sides of (2.1). Starting from $I$,
using (1.9), we get
\begin{align*}
I
&=
\prod_{k,j \geq 1}     
\left(
\frac{1}{1-z^{kj}} \right)^{{\mathcal C}(j,{\mathcal T}_{1}) {\mathcal C}(k,{\mathcal T}_{2})}
=\prod_{j \geq 1}     
\prod_{k \geq 1} 
\left(
\frac{1}{1-z^{jk}} 
\right)^{{\mathcal C}(k,{\mathcal T}_{2}){\mathcal C}(j,{\mathcal T}_{1})} =II.
\end{align*}
$\Box$

\noindent{\bf Remark 2.1.}  
Let ${\mathcal G}_{1}$ and  ${\mathcal G}_{2}$ be two graphs of the type of $G_{1}$ in Figure 1 with $R_{1}$ and $R_{2}$ oriented loops and directed edge adjacency matrices $T_{d1}$ and  $T_{d2}$, respectively. In this case,  (2.1) becomes the Strehl identity.

In the sequel
 we use ideas from [2] and [18] rather than the combinatorial arguments of Metropolis and Rota in [15] to prove several relations which are satisfied by ${\mathcal C}(N, {\mathcal T})$.

The relations we are about to prove involve the {\it Kronecker product} of adjacency matrices, hence, they  can be understood as relations for computing the number of cycles of a given length in the graph product in terms of the numbers of cycles in each graph.

Given two graphs
${\mathcal G}_{1}$ and  ${\mathcal G}_{2}$ with adjacency matrices ${\mathcal T}_{1}$ and  ${\mathcal T}_{2}$, respectively, the {\it Kronecker product} $ {\mathcal G}_{1} \otimes  {\mathcal G}_{2}$ is a graph with adjacency matrix ${\mathcal T}_{1} \otimes {\mathcal T}_{2}$. This is the matrix having the element $({\mathcal T}_{1})_{ij}$ replaced by the matrix $({\mathcal T}_{1})_{ij} {\mathcal T}_{2}$
. Many properties of the Kronecker product of graphs have been proved since Weichsel introduced it in [31] for non directed graphs. The product of directed graphs was investigated in [14]. 
A basic result is that the product graph need not be (strongly) connected even if graphs ${\mathcal G}_{1}$ and  ${\mathcal G}_{2}$ are. An explicit formula for the number of components is given in [31] and [14] together with conditions for  the product to be connected. The Ihara zeta function on Kronecker products have been investigated in [20].

Motivated by the structure of their relations Metropolis and Rota gave them  a ring theoretical interpretation which they investigated in connection with the {\it necklace ring}, {\it unital series} and {\it Witt vectors}. The relations below possibly will have a similar interpretation.

\begin{thm}
{\it Given matrices ${\mathcal T}_{1}$ and ${\mathcal T}_{2}$ and a positive integer $N$ define
\begin{equation}
S (N,{\mathcal T}_{i})=
 \sum_{g|N} \mu(g) \hspace{1mm} \Tr {\mathcal T}_{i}^{\frac{N}{g}}, i=1,2,
 \end{equation}
and 
denote by ${\mathcal T}_{1} \otimes {\mathcal T}_{2}$ the Kronecker product of ${\mathcal T}_{1}$ and ${\mathcal T}_{2}$. Then,
\begin{equation} 
\sum_{[s,t]=N} S (s,{\mathcal T}_{1}) S (t,{\mathcal T}_{2})
= S (N,{\mathcal T}_{1} \otimes {\mathcal T}_{2}).
\end{equation}
The summation is over the set of positive integers $\{s,t \mid [s,t]=N\}$,  $[s,t]$ being the least common multiple of $s$ and $t$. Furthermore, we have
\begin{equation} 
S (N,{\mathcal T}^{l})=\sum_{[l,t]=N l} S (t,{\mathcal T}).
\end{equation}  }
\end{thm}

\noindent{\bf Proof.} In order to prove  (2.3) it suffices to consider the equivalent formula (see [2])
\begin{equation*}
\sum_{k|N}\sum_{[s,t]=k } S (s,{\mathcal T}_{1}) S (t,{\mathcal T}_{2})
= \sum_{k|N} S (k,{\mathcal T}_{1} \otimes {\mathcal T}_{2}).
\end{equation*}
Using the M\"obius inversion formula, the left hand side is equal to
\begin{equation*}
\sum_{s|N} S (s,{\mathcal T}_{1}) \sum_{t|N} S (t,{\mathcal T}_{2})
= (Tr {\mathcal T}_{1}^{N}) (Tr {\mathcal T}_{2}^{N}).
\end{equation*}
But
$(\Tr {\mathcal T}_{1}^{N}) (\Tr {\mathcal T}_{2}^{N})=\Tr ({\mathcal T}_{1} \otimes {\mathcal T}_{2})^{N}$. By the
M\"obius inversion formula this gives the right hand side of the equivalent formula.
Using ideas from [18],
the next identity can be proved using the following equivalent formula:
\begin{equation*}
\sum_{g|N} S (N,{\mathcal T}^{l})
=\sum_{g|N} \sum_{[l,t]=\frac{Nl}{g}} S (t,{\mathcal T}).
\end{equation*}
The right hand side is equal to
$\sum_{t|lN} S(t,{\mathcal T})=\Tr {\mathcal T}^{lN}=\Tr ({\mathcal T}^{l})^{N}$.
Apply the M\"obius inversion formula to get the result.

\noindent{\bf Remark 2.2.} 
Formula (2.3) may be generalized to the case ${\mathcal T}={\mathcal T}_{1} \otimes {\mathcal T}_{2} \otimes ...\otimes {\mathcal T}_{l}$ to give
\begin{equation} 
\sum_{[{s}_{1}, \dots ,s_{l}]=N } S (s_{1},{\mathcal T}_{1}) \dots S (s_{l},{\mathcal T}_{l})
=  S (N,{\mathcal T}).
\end{equation}
Also, it can be proved that
\begin{equation} 
S(N,{\mathcal T}_{1}^{s} \otimes {\mathcal T}_{2}^{r})=
\sum_{[rp,sq]=Nrs} S (p,{\mathcal T}_{1})S (q,{\mathcal T}_{2}),
\end{equation}
where $r$ and $s$ are relatively prime and the summation is over all positive integers $p$ and $q$ such that $[rp,sq]=Nrs$. The proof is an application of previous identities as in [15], Theorem 5.

\begin{thm} 
 {\it Let $(s,t)$ denote the maximum common divisor of the positive integers $s$ and $t$. Then,
\begin{equation}
\sum_{[s,t]=N} (s,t) {\mathcal C}(s,{\mathcal T}_{1}) {\mathcal C} (t,{\mathcal T}_{2})
= {\mathcal C} (N,{\mathcal T}_{1} \otimes {\mathcal T}_{2}).
\end{equation}
The summation is over the set of positive integers $\{s,t \mid [s,t]=N\}$,  $[s,t]$ being the least common multiple of $s$ and $t$. Also,
\begin{equation} 
{\mathcal C} (N,{\mathcal T}^{l})=\sum_{[l,t]=N l} \frac{t}{N} {\mathcal C} (t,{\mathcal T}).
\end{equation}
and
\begin{equation} 
(r,s) {\mathcal C} \left(N,{\mathcal T}_{1}^{s/(r,s)} \otimes {\mathcal T}_{2}^{r/(r,s)} \right)=
 \sum (rp,sq) {\mathcal C} (p,{\mathcal T}_{1}) {\mathcal C}(q,{\mathcal T}_{2}).
\end{equation}
The sum is over $p$ and $q$  such that $pq/(pr,qs)=N/(r,s)$.}
\end{thm}

\noindent{\bf Proof.}
Use that $S(s,{\mathcal T})=s{\mathcal{C}}(s,{\mathcal T})$ and $[s,t](s,t)=st$ to get (2.7), and (2.8) also follows.
In terms of ${\mathcal C}$, 
(2.3) becomes
\begin{equation*}
N {\mathcal C} (N,{\mathcal T}_{1}^{s} \otimes {\mathcal T}_{2}^{r})=
\sum_{[rp,sq]=Nrs} pq {\mathcal C} (p,{\mathcal T}_{1}) {\mathcal C}(q,{\mathcal T}_{2}).
\end{equation*}
Using $(rp,sq)[rp,sq]=rpsq$ with $[rp,sq]=Nrs$ implies $(rp,sq)N=pq$ and
\begin{equation*}
{\mathcal C} (N,{\mathcal T}_{1}^{s} \otimes {\mathcal T}_{2}^{r})=
\sum_{[rp,sq]=Nrs} (rp,sq) {\mathcal C} (p,{\mathcal T}_{1}) {\mathcal C}(q,{\mathcal T}_{2}).
\end{equation*}
Replace $s$ and $r$ by $s/(r,s)$ and $r/(r,s)$ to get (2.9).

\noindent{\bf Remark 2.3.} 
Identity (2.7) can be extended to the general case ${\mathcal T}={\mathcal T}_{1} \otimes {\mathcal T}_{2} \otimes ...\otimes {\mathcal T}_{l}$ to give
\begin{equation} 
\sum_{[{s}_{1}, \dots ,{s}_{l}]=N } ({s}_{1}, \cdots, {s}_{l}  ){\mathcal C} ({s}_{1},{\mathcal T}_{1}) \dots {\mathcal C} ({s}_{l},{\mathcal T}_{l})
=  {\mathcal C} (N,{\mathcal T}),
\end{equation}
where $(s_{1},...,s_{l})$ is the greatest common divisor of $s_{1}$, ..., $s_{l}$ and the sum runs over all positive integers $s_{1}$, ..., $s_{l}$ with least common multiple $[s_{1}, \dots ,s_{l}]$ equal to $N$.

\noindent{\bf Remark 2.4.} 
 Let ${\mathcal G}_{1}$ and  ${\mathcal G}_{2}$ be two oriented graphs of the type of $G_{1}$ in Figure 1 with $R_{1}$ and $R_{2}$ oriented loops and directed adjacency matrices ${T}_{d1}$ and  ${T}_{d2}$, respectively. Then, ${\mathcal C} (N,T_{di})= \theta_{d}(N, T_{di})= {\mathcal M}(N;R_{i})$, $i=1,2$ and ${\mathcal C} (N,T_{d1} \otimes T_{d2})= {\mathcal M}(N;R_{1}R_{2})$. The above relations become the Metropolis-Rota identities. The graph product in this case has $R_{1}R_{2}$ edges hooked to a single vertex.

\section{Graph zeta functions and free Lie superalgebras} 

In this section we
connect the identity (1.9) and its reciprocal to exponentials, formal Taylor expansions and explicit formulas for the coefficients. These are given in Theorem 3.1 below. The identity (1.9) is included in the theorem for natural reasons. Theorem 3.1  is  important for the algebraic interpretation of  the relations (1.8) and (1.9), respectively, as a dimension formula and as the denominator (or generalized Witt) identity of a Lie superalgebra. 
In Theorem 3.2 we collect several recursions relating the coefficients of the Taylor expansions and the exponents ${\mathcal C}(N, {\mathcal T})$. Remarks 3.1-4 give additional informations.

Let ${\mathcal T}$, ${\mathcal C}(N,{\mathcal T})$, and $D$ be as  
in the Definition in the end of section 1.

\begin{thm} 
{\it Define
\begin{equation} 
 g(z):=\sum_{N=1}^{\infty} \frac{ \Tr {\mathcal T}^{N} }{N} z^{N}.
 \end{equation}
Then,
\begin{equation} 
\prod_{N=1}^{+\infty}  (1-z^{N})^{\pm {\mathcal C}(N,{\mathcal T})}
 = 
e^{\mp g(z)}  =[\det(1-z{\mathcal T})]^{\pm}= 1 \mp \sum_{i=1}^{+\infty} {d}_{\pm}(i) z^{i},\\
\end{equation}
 and 
\begin{equation} 
{d}_{\pm}(i)= \sum_{m=1}^{i} \lambda_{\pm}(m) \sum_{
\begin{array}{l} {a}_{1}+2{a}_{2}+...+i{a}_{i} =i\\
{a}_{1}+...+{a}_{i} = m \end{array}} 
 \prod_{k=1}^{i} 
\frac{(\Tr {\mathcal T}^{k})^{{a}_{k}}}{{a}_{k}! k^{{a}_{k}}},
\end{equation}
with $\lambda_{+}(m)=(-1)^{m+1}$, $\lambda_{-}(m)=+1$, $d_{+}(i)=0$ for $i > D$, and $d_{-}(i) \geq 0$, for all $i$'s. Furthermore,
\begin{equation}
\Tr {\mathcal T}^{N} = N  
\sum_{
\begin{array}{l}  s = ({s}_{i})_{i \geq 1}, {s}_{i} \in {\bf Z}_{\geq 0}\\
 \sum i{s}_{i}=N  \end{array}}
(\pm 1)^{|s|+1}
\frac{(\mid s \mid -1)!}{s!} \prod {d}_{\pm}(i)^{{s}_{i}},
\end{equation}
where
$\mid s \mid = \sum s_{i}, s! = \prod s_{i} !$.}
\end{thm}

\noindent{\bf Proof.}
Define $P_{\pm}$ by 
\begin{equation*}
{P}_{\pm}(z)=\prod_{N'=1}^{+\infty} (1-z^{N'})^{\pm {\mathcal C}(N',{\mathcal T})}.
\end{equation*}
Take the logarithm of both sides and use (1.8) to get
\begin{align*}
ln {P}_{\pm} &= \mp \sum_{N'} \sum_{k} \frac{1}{k} {\mathcal C}(N',{\mathcal T}) z^{N' k}
=\mp \sum_{N=1}^{+\infty} \sum_{k|N}\frac{1}{k} {\mathcal C}\left(\frac{N}{k}, {\mathcal T} \right) z^{N}\\
&= \mp \sum_{N=1}^{+\infty} \frac{\Tr {\mathcal T}^{N}}{N} z^{N} = \mp g(z)
\end{align*}
from which the first equality in (3.2) follows. 
From  the definition of $g(z)$, it follows that
\begin{align*}
\mp g(z):= \mp \sum_{N=1}^{\infty} \frac{\Tr {\mathcal T}^{N}}{N} z^{N} &= \mp \Tr \sum_{N=1}^{+\infty} \frac{1}{N} {\mathcal T}^{N}  z^{N}
= \pm \Tr \hspace{1mm}  ln(1-z{\mathcal T})\\
&= \pm ln \hspace{1mm}  \det(1-z{\mathcal T})\\
\end{align*}
proving the second equality in (3.2).

The third equality in (3.2) is obtained by formally expanding the exponential.
As the formal Taylor
expansion of $1-e^{\mp g}$, the coefficients  $c_{\pm}$ are given by
\begin{equation*}
{d}_{\pm}(i) =
\frac{1}{i!}
\frac{d^{i}}{d z^{i}} \left[ \pm(1-e^{\mp g}) \right]|_{z=0}.
\end{equation*}
Using Faa di Bruno's formula as in [4,5],  the derivatives  can be computed explicitly and (3.3) follows. The determinant is a polynomial of maximum degree $D$, hence, $d_{+}(i)=0$ for $i>D$. Clearly, $d_{-}(i) \geq 0$.

To prove (3.4) write 
\begin{equation*}
{\mathcal F}:=\mp ln \left( 1 \mp \sum_{i} {d}_{\pm}(i) z^{i} \right) 
= \mp \sum_{l=1}^{+\infty} \frac{1}{l} \left( \pm \sum_{i} {d}_{\pm}(i) z^{i} \right)^{l}.
\end{equation*}
Expand the right hand side in powers of $z$ to get:
\begin{align*}
 {\mathcal F} &= \pm \sum_{l=1}^{+\infty} \frac{(\pm 1)^{l}}{l} 
 \sum_{
\begin{array}{l}  s = ({s}_{i})_{i \geq 1}, {s}_{i} \in {\bf Z}_{\geq 0}\\
 \sum i{s}_{i}=l  \end{array}}
 \frac{(\sum {s}_{i})!}{\prod ({s}_{i} !)} \left(\prod {d}_{\pm}(i)^{{s}_{i}} \right) z^{\sum {s}_{i} l}\\
 &= \sum_{k=1}^{+\infty}    z^{k} 
 \sum_{
\begin{array}{l}  s = ({s}_{i})_{i \geq 1}, {s}_{i} \in {\bf Z}_{\geq 0}\\
 \sum i{s}_{i}=k  \end{array}}
 (\pm 1)^{|s|+1} \frac{(\mid s \mid -1)!}{s!} \prod {d}_{\pm}(i)^{{s}_{i}}.
\end{align*}
The second equality in (3.2) applied to the left hand side yields
\begin{equation*}
{\mathcal F} =\sum_{k=1}^{+\infty} \frac{\Tr {\mathcal T}^{k}}{k} z^{k}.
\end{equation*}
Comparing the coefficients give the result.

The general expressions (3.3) for the coefficients $d_{\pm}(i)$ given in  Theorem 3.1 are complicated. They can be computed recursively as the
next theorem shows.

\begin{thm} {\it Set $\omega(n):= \Tr {\mathcal T}^{n}$. Then,
\begin{align} 
{d}_{\pm} (1) &= \omega(1), \\
 n {d}_{\pm} (n) &=  \omega(n) \mp \sum_{k=1}^{n-1}  \omega(n-k) {d}_{\pm}(k), n \geq 1,\\
 {d}_{-}(n) &= {d}_{+}(n)+\sum_{i=1}^{n-1} {d}_{+}(i) {d}_{-}(n-i), n \geq 1, \\
 |{d}_{+}(n)| & \leq   {d}_{-}(n),\\
 {\mathcal C}(n, {\mathcal T}) &= {d}_{+}(n)+\frac{1}{n} \sum_{k=1}^{n-1}\left( \sum_{g \mid k}
 g {\mathcal C}(g, {\mathcal T})\right) {d}_{+}(n-k) 
 - \sum_{n \neq g \mid n}
 \frac{g}{n} {\mathcal C}(g, {\mathcal T}) .
\end{align}  
}
\end{thm}

\noindent{\bf Proof.} Define $f(z):=e^{\mp h(z)}$, $g_{n} := \mp \omega_{n}$, and denote by  $f'$ the formal derivative of $f$. Also, define $c(n):=\mp d_{\pm}(n)$.
Then,
\begin{equation*} 
f'(z)= f(z) \sum_{n=1}^{\infty} g(n) z^{n-1} = \sum_{n=1}^{\infty} c(n) n z^{n-1} .
\end{equation*}
Thus,
\begin{equation*}
\left( 1+\sum_{n=1}^{\infty} c(n) z^{n} \right) \left( \sum_{n=1}^{\infty} g(n) z^{n-1} \right)= \sum_{n=1}^{\infty} c(n) n z^{n-1}.
\end{equation*}
Equating the coefficients of both sides, we get
\begin{equation*} 
n c(n)= g(n)+ \sum_{k=1}^{n-1} g(n-k) c(k).
\end{equation*}
From the definitions the result follows. Now, use that $e^{-h(z)} e^{+h(z)}=1$, to prove (3.7).

To prove (3.9) add $d_{+}(n)$ to $d_{-}(n)$ to get that
\begin{equation*}
{d}_{-}(n)+{d}_{+}(n)=\frac{w(n)}{n} \geq 0.
\end{equation*}
Then, using this result, subtract $d_{+}$ from $d_{-}$  to get
\begin{equation*}
n(d_{-}(n)-d_{+}(n))= \sum_{k=1}^{n-1} w(n-k)(d_{-}(k)+d_{+}(k)) \geq 0.
\end{equation*}
Relation (3.9) follows from  (3.6) using that $\omega(n)=\sum_{g\mid n} g {\mathcal C}(g, {\mathcal T})$.

In two papers S. -J. Kang and M. -H Kim [7,8] generalized the
Proposition 1.1 (section 1)
to the case that  the free Lie algebra $L$ is generated by an infinite graded vector space.  They obtained  a generalized Witt formula for the dimensions of the homogeneous subspaces of $L$ which satisfies a  generalized Witt identity that plays the role of the denominator identity for the free Lie algebra.
 In [9] S. -J. Kang extended the results to superspaces and Lie superalgebras. 
The following proposition summarizes Kang's results from [9] which are relevant for our objectives.

\begin{prop} 
{\it  Let ${\mathcal V}= \bigoplus_{N=1}^{\infty}
{\mathcal V}_{N}$ be a ${\mathbb{Z}}_{>0}$-graded superspace  with finite dimensions $\dim {\mathcal V}_{N}= |t_{N}|$ and superdimensions
$\Dim {\mathcal V}_{N}= t_{N} \in {\mathbb{Z}}$, $\forall N \geq 1$. 
Let ${\mathcal L}=
\bigoplus_{N=1}^{\infty} {\mathcal L}_{N}$ be the free Lie superalgebra generated
by ${\mathcal V}$ with a ${\mathbb{Z}}_{>0}$-gradation induced by that of ${\mathcal V}$. Then, the ${\mathcal L}_{N}$
superdimension  is
\begin{equation} 
\Dim {\mathcal L}_{N}= \sum_{g | N} \frac{\mu(g)}{g}  W \left(\frac{N}{g}\right).
\end{equation}
The summation ranges over all positive divisors $g$ of $N$ and $W$ is given by
\begin{equation} 
 W(N)= \sum_{s \in T(N)} \frac{(\mid s \mid -1)!}{s!} \prod t(i)^{{s}_{i}},
\end{equation}
where $T(N)=\{ s = (s_{i})_{i \geq 1} \mid s_{i} \in {\mathbb{Z}}_{\geq 0}, 
\sum is_{i}=N \}$
and $\mid s \mid = \sum s_{i}$, $ s! = \prod s_{i} !$. Furthermore,
\begin{equation}
\prod_{N=1}^{\infty} (1-z^{N})^{\pm \Dim {\mathcal L}_{N}}= 1 \mp \sum_{N=1}^{\infty} {f}_{\pm}(N) z^{N},
\end{equation}
with $f_{+}(N)=t_{N}$ and $f_{-}(N)=\Dim U({\mathcal L})_{N}$,
where $\Dim U({\mathcal L})_{N}$ is the dimension of the $N$-th homogeneous subspace of the
 universal enveloping algebra $U({\mathcal L})$ and the generating function for the $W$'s,
\begin{equation}
g(z) :=\sum_{N=1}^{\infty} W(N)z^{N},
\end{equation}
satisfies
\begin{equation}
e^{-g(z)}=  1-\sum_{N=1}^{\infty} {t}_{N}z^{N}.
\end{equation}  }
\end{prop}

\noindent{\bf Remark 3.1.} 
 In [9], (3.10) is called  the {\it generalized Witt formula}; $W$ is
called the Witt partition function; and the $(+,-,+)$ case of (3.12) is called the denominator identity (or the generalized Witt identity, in the present paper) of the free Lie superalgebra.

In the sequel we apply ideas from section  2.3 of [9] to interpret algebraically the relations (1.8-9).

Given a formal power series $\sum_{N=1}^{+\infty} t_{N} z^{N}$ with $ t_{N} \in {\mathbb Z}$, for all $i \geq 1$,
the coefficients in the series can be interpreted as the superdimensions of a ${\mathbb Z}_{>0}$-graded  superspace ${\mathcal V}= \bigoplus_{i=N}^{\infty}
{\mathcal V}_{N}$ with dimensions $dim {\mathcal V}_{N}= |t_{N}|$ and superdimensions
$\Dim {\mathcal V}_{N}= t_{N} \in {\mathbb Z}$.
Let ${\mathcal L}$ be the free Lie superalgebra generated by ${\mathcal V}$. Then, it has a gradation induced by ${\mathcal V}$ and its homogeneous subspaces have dimension given by (3.10) and (3.11).  Apply this interpretation to the  determinant $\det(1-z{\mathcal T})$ which is a polynomial in the formal variable $z$ of degree $D$. Recall that $D=2|E|$, if $\mathcal T=T$; $D=|V|$, if $\mathcal T=A_{d}$; and $D=|E|$, if $T=T_{d}$. 
The polynomial is a power series in the variable $z$ with coefficients $t_{N}=0$, for $N > D$.
Comparison of the formulas
in Theorem 3.1 with the formulas
in the above Proposition yields that
given a graph $G$, ${\mathcal T}$ its associated matrix, 
let
${\mathcal V}= \bigoplus_{N=1}^{D} 
{\mathcal V}_{N}$ be a ${\mathbb{Z}}_{>0}$-graded superspace with finite dimensions 
$\dim {\mathcal V}_{N}= |d_{+}(N)|$ and the superdimensions
$\Dim {\mathcal V}_{N}= d_{+}(N)$ given by (3.3), the coefficients of $\det(1-z{\mathcal T})$. Let ${\mathcal L}=
\bigoplus_{N=1}^{\infty} {\mathcal L}_{N}$ be the free Lie superalgebra generated
by ${\mathcal V}$. Then, the  ${\mathcal L}_{N}$
superdimension is  $\Dim{\mathcal L}_{N} ={\mathcal C}(N, T)$ given by (1.8). This means that the equivalence classes of nonperiodic cycles of length $N$ in $G$ form a basis
of ${\mathcal L}_{N}$.
The algebra has denominator (or generalized Witt) identity given by (1.9). 

The following new interpretation of 
the Ih or BL zeta functions of a graph follows. The zeta function  $\zeta_{G}(z) := det(1-z{\mathcal T})^{-1}$
of a graph $G$ is the generating function for the dimensions $\Dim U({\mathcal L})_{N}=d_{-}(N)$, given by  (3.3),  of the subspaces of
the enveloping algebra $U({\mathcal L})$ of Lie superalgebra ${\mathcal L}$ generated by G.

\noindent{\bf Remark 3.2} 
In [10], S.-J. Kang, J.-H. Kwon, and Y.-T. Oh derived Peterson-type dimension formulas for graded Lie superalgebras. In particular, see the Example 3.6, p. 118 of [10]. Then, the Peterson-type formulas are recursive relations between the dimensions
$\Dim {\mathcal L}_{N}$ and the coefficients $f_{+}(N)$ in (3.12). See the formula after relation (3.15) in [10]. Using  our notation, this formula is exactly the relation (3.9) in  Theorem 3.2.

\noindent{\bf Remark 3.3}
Given two graphs and the free Lie superalgebras generated by them the algebra
generated by the Kronecker product graph
will have dimensions that can be expressed in terms of the dimensions of the algebras generated by  the individual graphs. They are given by the Carlitz-Metropolis-Rota-type identities derived in section 2. It would be interesting to investigate how the vector spaces that generate the algebras of the individual graphs are related to the vector space that generate the algebra of the graph product.

\noindent{\bf Remark 3.4} 
It is known from the works of several authors that two graphs can have the same Ihara zeta function. See [24-25] and references therein. This is also true for the Bowen-Lanford zeta function. See Example 4.1.
Graphs with the same zeta function will generate the same algebra. Since the functions have the same coefficients $d_{+}$ it follows from (3.9) that they  have the same numbers of cycles of same length. 
This is the case of the graphs $G_{1}$ in Figure 1 (section 1) with $R=2$ and the graph $G_{2}$ in Figure 2 below.
The  edge and vertex adjacency matrices of $G_{2}$ are
\begin{equation*}
{T}_{d} = \left( \begin{array}{clcrcl}
1 & 1 & 0 & 0 & 0 & 0 \\
0 & 0 & 1 & 1 & 0 & 0 \\
0 & 0 & 1 & 1 & 0 & 0 \\
0 & 0 & 0 & 0 & 1 & 1 \\
0 & 0 & 1 & 1 & 0 & 0 \\
1 & 1 & 0 & 0 & 0 & 0 \\
\end{array} \right),
\hspace{5mm} 
{A}_{d}= \left( \begin{array}{clcrcl}
0 & 1 & 1 \\
1 & 1 & 0 \\
0 & 1 & 1  \\
\end{array} \right).
\end{equation*}
The loops $1$ and $3$ are hooked to  vertices $3$ and $2$, respectively. The matrices have the traces
$\Tr T_{d}^{N}=\Tr A_{d}^{N}=
2^{N}$, and
\begin{equation*}
\det(1-z {T}_{d})= \det(1-z {A}_{d}) = 1-2 z.
\end{equation*}
\begin{center}
\begin{figure}[ht]
\centering
\includegraphics[scale=0.5]{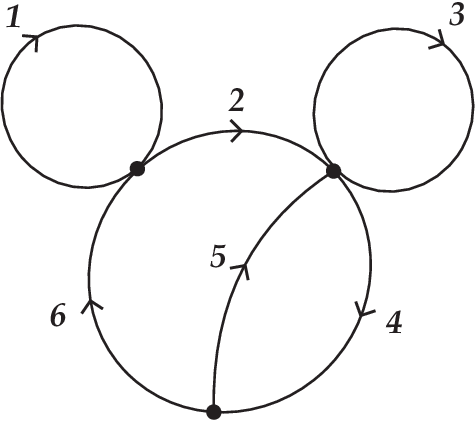}
\caption{Graph $G_{2}$}
\label{Fi:G2}
\end{figure}
\end{center}
The number of classes of nonperiodic cycles of length $N$ is
\begin{equation*}
\theta_{d} (N)= \frac{1}{N}\sum_{g|N} \mu (g) 2^{\frac{N}{g}}.
\end{equation*}
and
\begin{equation*}
\prod_{N=1}^{+\infty} (1-z^{N})^{\theta_{d}(N)}
=  1-2 z.
\end{equation*}
where $\theta_{d}(N)=\theta_{d}(N,T_{d})=\theta_{d}(N,A_{d})$.

\section{Examples}

\noindent{\bf Example 4.1.} 
$G_{1}$, the graph with $R \geq 2$ edges hooked to a single vertex is shown in Figure 1, section 1. Let's orient the edges counterclockwisely.
Add in the inverse edges to get the edge adjacency matrix $T$ for the non oriented  $G_{1}$. This is the $2R \times 2R$ symmetric matrix
\begin{equation*}
{T} = \left( \begin{array}{clcr}
A & B\\
B & A 
\end{array} \right)
\end{equation*}
where $A$ is the $R \times R$ matrix with all entries equal to $1$ and $B$  is the $R \times R$ matrix with  the main diagonal entries  equal to $0$ and all the other entries equal to $1$. This matrix has trace given by
\begin{equation*}
\Tr {T}^{N} = 1+(R-1)(1+(-1)^{N})+(2R-1)^{N}, \hspace{2mm} N=1,2, \dots,
\end{equation*}
and determinant
\begin{eqnarray*}
\det(1-z{T}) &=& (1-z) \left[ 1-(2R-1)z \right](1-z^{2})^{R-1}\\
            &=& 1-\sum_{i=1}^{2R} c(i) z^{i},
\end{eqnarray*}
where $c(2R)=(-1)^{R}(2R-1)$,
\begin{equation*}
c(2i)=(-1)^{i}(2i-1) 
\left( \begin{array}{c}
R\\
i
\end{array} \right), \hspace{2mm} i=1, \cdots, R-1,
\end{equation*}
and
\begin{equation*}
c(2i+1)=2R(-1)^{i} \left( \begin{array}{c}
R-1\\
i
\end{array} \right), \hspace{2mm} i=0,1, \cdots R-1.
\end{equation*}
Futhermore,
\begin{eqnarray*}
[det(1-z{T})]^{-1} &=& \sum_{q=0}^{+\infty} z^{q} \sum_{i=0}^{q} {a}_{i} (2R-1)^{q-i},
\end{eqnarray*}
where
\begin{equation*}
{a}_{i}= \sum_{k=0}^{i} (-1)^{i-k} 
\left( \begin{array}{c}
k+R-1\\
R-1
\end{array} \right)
\left( \begin{array}{c}
i-k+R-2\\
R-2
\end{array} \right).
\end{equation*}
Let's consider the case $R=2$. In this case,
\begin{equation*}
\Tr {T}^{N}=2+(-1)^{N}+3^N, \hspace{5mm} \det(1-z{T}_{{G}_{1}})= 1-4z+2z^2+4z^3-3z^4,
\end{equation*}
so that the number of classes of nonperiodic cycles of length $N$ is given by the formula
\begin{equation*}
\theta(N, {T}) = \frac{1}{N} \sum_{g|N} \mu (g) \left( 2+(-1)^{\frac{N}{g}}+3^{\frac{N}{g}} \right).
\end{equation*}
The graph generates the following algebra.
Let
${\mathcal V}= \bigoplus_{i=1}^{4}
{\mathcal V}_{i}$ be a ${\mathbb Z}_{>0}$-graded superspace  with dimensions  $\dim {\mathcal V}_{1}=4$, $\dim {\mathcal V}_{2}=2$, $\dim {\mathcal V}_{3}=4$, $\dim {\mathcal V}_{4}=3$ and superdimensions $\Dim {\mathcal V}_{1}=-4$, $\Dim {\mathcal V}_{2}=2$, $\Dim {\mathcal V}_{3}=4$, $\Dim {\mathcal V}_{4}=-3$.  Let ${\mathcal L}=
\bigoplus_{N=1}^{\infty} {\mathcal L}_{N}$ be the free graded Lie superalgebra generated
by ${\mathcal V}$.  The
dimensions of the subspaces are
\begin{equation*}
\Dim{\mathcal L}_{N} = \theta(N,T), N=1,2,....
\end{equation*}
For instance, $\Dim {\mathcal L}_{1}=4$, $\Dim {\mathcal L}_{2}=4$, $\Dim {\mathcal L}_{3}=8$, 
with basis:
\begin{eqnarray*}
{\mathcal L}_{1}  & : & \{[{e}_{1}], [{e}_{1}^{-1}], [{e}_{2}], [{e}_{2}^{-1}] \}, \\
{\mathcal L}_{2}  & : & \{[{e}_{1}{e}_{2}], [{e}_{1}{e}_{2}^{-1}], [{e}_{1}^{-1}{e}_{2}], [{e}_{1}^{-1}{e}_{2}^{-1}] \},\\
 {\mathcal L}_{3}  & : &
\{[{e}_{1}^{+2}{e}_{2}^{+1}], [{e}_{1}^{+2}{e}_{2}^{-1}], [{e}_{1}^{-2}{e}_{2}^{+1}], [{e}_{1}^{-2}{e}_{2}^{-1}],  [{e}_{1}^{+1}{e}_{2}^{+2}], [{e}_{1}^{+1}{e}_{2}^{-2}], [{e}_{1}^{-1}{e}_{2}^{+2}], [{e}_{1}^{-1}{e}_{2}^{-2}]\}.
\end{eqnarray*}

The generalized Witt identity is
\begin{eqnarray*}
\prod_{N=1}^{+\infty} (1-z^{N})^{\theta(N, {T})}
&=&  1-4z+2z^2+4z^3-3z^4.
\end{eqnarray*}
 The dimensions of 
the subspaces $U_{n}({\mathcal L})$ of the enveloping algebra $U({\mathcal L})$ have dimensions generated by
the Ihara zeta function of the graph,
\begin{eqnarray*}
\prod_{N=1}^{+\infty} (1-z^{N})^{-\theta(N, {T}_{{G}_{1}})}
&=& 1+\frac{1}{16}\sum_{n=1}^{\infty} ((-1)^{n}+ 3^{n+3}-12-4n) z^{n}.
\end{eqnarray*}
The first few terms give $\Dim U_{1}=4$, $\Dim U_{2}=14$, $\Dim U_{3}=44$, $\Dim U_{4}=135$.

\noindent{\bf Example 4.2.} 
$G_{3}$, the bipartite graph shown in Figure 3.
The edge adjacency matrix of $G_{3}$ is
\begin{equation*}
{T} = \left( \begin{array}{clcrclcr}
0 & 1 & 0 & 0 & 0 & 1\\
1 & 0 & 1 & 0 & 0 & 0 \\
0 & 1 & 0 & 1 & 0 & 0 \\
0 & 0 & 1 & 0 & 1 & 0\\
0 & 0 & 0 & 1 & 0 & 1\\
1 & 0 & 0 & 0 & 1 & 0 
\end{array} \right).
\end{equation*}
The matrix has the trace
$\Tr T^{N}= 0$ if $N$ is odd and $\Tr T^{N}=
4+2 \cdot 2^{N}$ if $N$ is even, and the determinant
\begin{equation*}
\det(1-z{T})= 1-(6z^2-9z^4+4z^6).
\end{equation*}
\begin{center}
\begin{figure}[ht]
\centering
\includegraphics[scale=0.5]{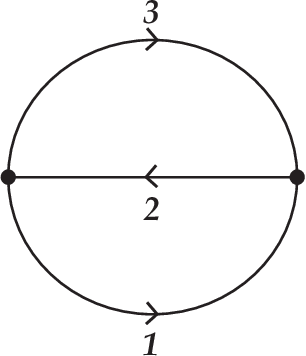}
\caption{Graph $G_{3}$}
\label{Fi:G3}
\end{figure}
\end{center}
The number of classes of nonperiodic cycles of length $N$ is $\theta(N, T )=0$, if $N$ is odd, and
\begin{equation*}
\theta (N, {T})= \frac{1}{N}
\sum_{\substack{g|N \\ N/g \hspace{1mm} even}} \mu (g) \left( 4+2^{\frac{N}{g}+1}
\right),
\end{equation*}
if $N$ is even. The graph generates the following algebra.
Let
${\mathcal V}= \bigoplus_{i=1}^{3}
{\mathcal V}_{2i}$ be a ${\mathbb Z}_{>0}$-graded superspace  with dimensions  $\dim {\mathcal V}_{2}=6$, $\dim {\mathcal V}_{4}=9$, $\dim {\mathcal V}_{6}=4$ and superdimensions $\Dim {\mathcal V}_{2}=6$, $\Dim {\mathcal V}_{4}=-9$, $\Dim {\mathcal V}_{6}=4$.  Let ${\mathcal L}=
\bigoplus_{N=1}^{\infty} {\mathcal L}_{N}$ be the free graded Lie superalgebra generated
by ${\mathcal V}$.  The
dimension of ${\mathcal L}_{N}$ subspace is $\Dim{\mathcal L}_{N}=\theta(N, T)$.
For instance, $\Dim {\mathcal L}_{2}=6$ and $\Dim {\mathcal L}_{4}=6$. The
 basis are:
\begin{eqnarray*}
{\mathcal L}_{2}  & : & \{
[{e}_{1} {e}_{2}], [{e}_{1}^{-1} {e}_{2}^{-1}], [{e}_{1} {e}_{3}^{-1}], [{e}_{1}^{-1} {e}_{3}], [{e}_{2} {e}_{3}],  [{e}_{2}^{-1} {e}_{3}^{-1}] \},\\
 {\mathcal L}_{4}  & : &
\{
[{e}_{1} {e}_{2} {e}_{3} {e}_{2}],
[{e}_{1} {e}_{2} {e}_{1} {e}_{3}^{-1}], [{e}_{1} {e}_{3}^{-1} {e}_{2}^{-1} {e}_{3}^{-1}], [{e}_{1}^{-1} {e}_{2}^{-1} {e}_{3}^{-1}{e}_{2}^{-1}],
 [{e}_{1}^{-1} {e}_{2}^{-1} {e}_{1}^{-1}{e}_{3}], \\
 & &[{e}_{1}^{-1} {e}_{3} {e}_{2}{e}_{3}]
\}.
\end{eqnarray*}
The dimensions satisfy 
the generalized Witt identity
\begin{equation*}
\prod_{N=1}^{+\infty} (1-z^{N})^{\theta(N, {T})}
=  1-6z^2+9z^4-4z^6.
\end{equation*}
The generating function for the dimensions of the subspaces $U_{n}({\mathcal L})$ of the enveloping algebra $U({\mathcal L})$ is given by the Ihara zeta function of $G_{3}$:
\begin{align*}
\prod_{N=1}^{+\infty} (1-z^{N})^{-\theta(N, T)}
&= 1+\frac{1}{18}\sum_{n=1}^{\infty} (2^{2n+5}-6n-14) z^{2n}.
\end{align*}

\noindent{\bf Example 4.3.}
$G_{3}$, the  graph shown in Figure 3, is strongly connected. Call $v_{1}$ ($v_{2}$) the vertex on the left (right).
The directed edge and vertex adjacency matrices of $G_{3}$ are
\begin{equation*}
{{T}_{d}} = \left( \begin{array}{clcr}
0 & 1 & 0 \\
1 & 0 & 1 \\
0 & 1 & 0 \\
\end{array} \right),
\hspace{5mm}
{A}_{d} = 
\left( \begin{array}{clcr}
0 & 2 \\
1 & 0 \\
\end{array} \right).
\end{equation*}
The matrices have the traces
$\Tr { T_{d}}^{N} = \Tr A_{d}^{N}=0$ if $N$ is odd and $\Tr {T_{d}}^{N}=\Tr A_{d}^{N}=
2^{\frac{N}{2}+1}$ if $N$ is even, and the determinants
\begin{equation*}
\det(1-z { {T}_{d}})= \det(1-z {A}_{d}) = 1-2 z^2.
\end{equation*}
The number of classes of nonperiodic cycles of length $N$ is $\theta_{d}(N )=0$, if $N$ is odd, and
\begin{equation*}
\theta_{d} (N, T_{d})= \frac{1}{N}\sum_{\substack{g|N \\ N/g \hspace{1mm} even}} \mu (g)  2^{\frac{N}{2g}+1} ,
\end{equation*}
if $N$ is even. The first few values are $\theta_{d}(2)=2$, $\theta_{d}(4)=1$,
$\theta_{d}(6)=2$, $\theta_{d}(8)=3$, $\theta_{d}(10)=6$. 
For $N=2$, the classes are $[e_{1} e_{2}]$ and $[e_{2} e_{3}]$. For $N=4$,
only $[e_{1}e_{2}e_{3}e_{2}]$.

Let ${\mathcal V}$ be the vector space with $\dim {\mathcal V}=2$. It generates the graded free Lie algebra $\bigoplus_{N=1}^{\infty} {\mathcal L}_{N}$ with $\dim {\mathcal L}_{N}=0$, if $N$ is odd, and $\dim {\mathcal L}_{N}= \theta_{d}(N)$, if $N$ is even. For instance, $\dim {\mathcal L}_{2}=2$, $\dim {\mathcal L}_{4}=1$.
${\mathcal L}_{2}$ and ${\mathcal L}_{4}$ have basis
\begin{eqnarray*}
{\mathcal L}_{2}  & : & \{
[{e}_{1} {e}_{2}], [{e}_{2}{e}_{3}] \},\\
 {\mathcal L}_{4}  & : &
\{
[{e}_{1} {e}_{2} {e}_{3} {e}_{2}]
\}.
\end{eqnarray*}
The algebra has generalized Witt identity
\begin{equation*}
\prod_{N=1}^{+\infty} (1-z^{N})^{\theta_{d}(N)}
=  1-2z^2.
\end{equation*}
and the dimensions of the spaces $U_{n}(L)$ of the enveloping algebra are generated by the Bowen-Lanford zeta function
\begin{equation*}
\prod_{N=1}^{+\infty} (1-z^{N})^{-\theta_{d}(N)}
=  (1-2z^2)^{-1}=\sum_{n=0}^{\infty} 2^{n} z^{2n}.
\end{equation*}

\noindent{\bf Example 4.4.}
$G_{4}$, the  graph shown in Figure 4. Call $v_{4}$ the upper vertex and the others,  $v_{1},v_{2}, v_{3}$,  from left to right.
The directed edge adjacency  and the vertex adjacency matrices of $G_{4}$ are
\begin{equation*}
{{T}_{d}} = \left( \begin{array}{clcrcl}
0 & 1 & 1 & 0 & 0 & 0 \\
1 & 0 & 0 & 0 & 1 & 0 \\
0 & 0 & 0 & 1 & 0 & 0 \\
0 & 1 & 1 & 0 & 0 & 0 \\
0 & 0 & 0 & 0 & 0 & 1 \\
0 & 0 & 0 & 1 & 0 & 0 \\
\end{array} \right),
\hspace{5mm}
{A}_{d} = \left( \begin{array}{clcrcl}
0 & 1 & 0 & 1 \\
1 & 0 & 1 & 0 \\
0 & 1 & 0 & 0 \\
0 & 0 & 1 & 0 \\
\end{array} \right).
\end{equation*}
The matrices have  determinants
\begin{equation*}
\det{(1-z{T}_{d})}=\det{(1-z{A}_{d})}= 1-(2z^{2}+z^{4}).
\end{equation*}
If $N$ is odd, $\Tr {T_{d}}^{N}=0$. If $N$ is even, it follows from (4.4) that
\begin{equation*}
\Tr {{T}_{d}}^{N}=\Tr {{A}_{d}}^N = N  \sum_{2{a}+4{b}=N}
\frac{(a+b -1)!}{a!b!} {2}^{a}.
\end{equation*}
\begin{center}
\begin{figure}[ht]
\centering
\includegraphics[scale=0.5]{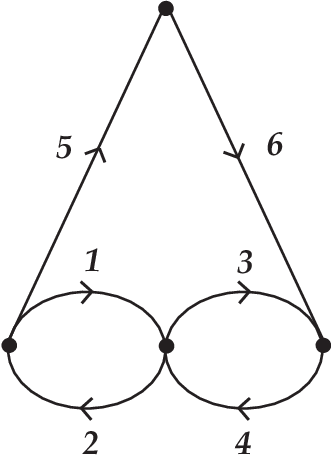}
\caption{Graph $G_{4}$}
\label{Fi:G4}
\end{figure}
\end{center}
The number of classes of nonperiodic cycles  for $N$ odd is zero. For $N$ even,
\begin{equation*}
\theta_{d} (N)= \sum_{\substack{g|N \\ N/g \hspace{1mm} even}}  \frac{\mu (g)}{g}  \sum_{2{a}+4{b}=\frac{N}{g}}
\frac{(a+b -1)!}{a!b!} {2}^{a}.
\end{equation*}
For instance, $\theta_{d}(2)=2$, $\theta_{d}(4)=2$,
$\theta_{d}(6)=2$. 
The cycles are $[e_{1}e_{2}]$ and $[e_{3}e_{4}]$, for $N=2$, $[e_{1}e_{3}e_{4}e_{2}]$ and $[e_{2}e_{5}e_{6}e_{4}]$, for $N=4$, and
$[e_{2}e_{5}e_{6}e_{4}e_{2}e_{5}]$, $[e_{2}e_{5}e_{6}e_{4}e_{3}e_{4}]$,
$[e_{1}e_{3}e_{4}e_{2}e_{1}e_{2}]$ and $[e_{1}e_{3}e_{4}e_{3}e_{4}e_{2}]$, for $N=6$.
Also,
\begin{equation*}
\prod_{N=1}^{+\infty} (1-z^{N})^{\theta_{d}(N)}
=  1-2z^{2}-z^{4}
\end{equation*}

Let
${\mathcal V}= \bigoplus_{i=1}^{4}$, $\dim {\mathcal V}_{1}=\dim {\mathcal V}_{3}=0$
$\dim V_{2}=2$, $\dim V_{4}=1$.
It generates the graded free Lie algebra $\bigoplus_{N=1}^{\infty} {\mathcal L}_{N}$ with $\dim {\mathcal L}_{N}=0$, if $N$ is odd, and $\dim {\mathcal L}_{N}= \theta(N)$, if $N$ is even. For instance, $\dim {\mathcal L}_{2}=2$, $\dim {\mathcal L}_{4}=2$.
${\mathcal L}_{2}$ and ${\mathcal L}_{4}$ have basis
\begin{eqnarray*}
{\mathcal L}_{2}  & : & \{
[{e}_{1} {e}_{2}], [{e}_{3}{e}_{4}] \},\\
 {\mathcal L}_{4}  & : &
\{
[{e}_{1} {e}_{3} {e}_{4} {e}_{2}], [{e}_{2} {e}_{5} {e}_{6} {e}_{4}]
\}.
\end{eqnarray*}
The algebra has generalized Witt identity
\begin{equation*}
\prod_{N=1}^{+\infty} (1-z^{N})^{\theta_{d}(N)}
=   1-2z^{2}-z^{4},
\end{equation*}
and the dimensions of the spaces $U_{n}(L)$ of the enveloping algebra are generated by the Bowen-Lanford zeta function
\begin{equation*}
\prod_{N=1}^{+\infty} (1-z^{N})^{-\theta_{d}(N)}
=  (1-2z^2-z^{4})^{-1}=\sum_{n=0}^{\infty} {d}_{-}(n) z^{n}.
\end{equation*}
The coefficients can be computed recursively in the following manner.
We have that $d_{+}(1)=0$, $d_{+}(2)=2$, $d_{+}(3)=0$, $d_{+}(4)=1$, and $d_{+}(n)=0$, if $n \geq 5$. Using the relation (4.8), theorem 4.2, we find that
\begin{equation*}
{d}_{-}(n)={d}_{+}(n)+{d}_{+}(2){d}_{-}(n-2)+{d}_{+}(4){d}_{-}(n-4)
\end{equation*}
from which we get that $d_{-}(n)=0$, if $n$ is odd and  ${d}_{-}(2)=2$, ${d}_{-}(4)=5$, ${d}_{-}(6)=12$, etc.

\section{Restricted necklace colorings} \label{sec:col}

In this section we interpret the cycle counting formula (1.8)  as a counting formula for the number of classes of non periodic colorings
of a necklace with $N$ beads. 

First, let's consider (1.8) in the case ${\mathcal T}=T$. Given a graph $G$ with $|E|$ edges and the colors $c_{1}$, ... , $c_{2|E|}$,
assign  $c_{i}$, $c_{|E|+i}$ to the edges $e_{i}$, $e_{|E|+i}=e_{i}^{-1} \in G'$, respectively, so that
to a cycle of length $N$ in $G$ corresponds with an ordered sequence of $N$ colors. Assign each color in this sequence to a bead in a circular string with $N$ beads - a necklace - in such a manner that two adjacent colors in the sequence are assigned to adjacent beads. The non backtracking condition for cycles implies that 
no two adjacent  beads are painted with colors, say, $c_i$ and $c_{|E|+i}$.
There is a 
 correspondence between the classes of nonperiodic cycles of length $N$ in $G$ and classes of  nonperiodic colorings of a necklace with $N$ beads with at most $2|E|$ distinct colors  induced by the cycles so that the number of inequivalent colorings is $\theta(N, T)$, given by (1.1). Of course,
the structure of the graph reflects itself  in the coloring so that the coloring is restricted by that structure.
For instance, the presence of loops in the graph means that their assigned colors  may appear repeated in a string of adjacent beads.  This can not happen to a color assigned to an edge which is not a loop. 
 The edge adjacency matrix $T$ may be called the {\it color matrix}. It basically tells what colors are allowed to be adjacent to a given color in the necklace. 
Element $T_{ij}=1$, if a color $c_j$ can be adjacent to color $c_i$ and $c_j \neq
c_{|E|+i}$;
$T_{ij}=0$, otherwise.

Inversely, any matrix $T$ with zeros and ones as entries
and even order and that has the correct structure (see [23], Lemma 4 on p. 151, and [6], Remark 1.5 on p. 7) is a necklace color matrix that can be interpreted as the edge adjacency matrix of some graph $G$, hence, $\theta (N, T)$ is the number of classes of  nonperiodic colorings in the necklace and of cycles in the graph.

The same ideas apply to  directed graphs. 
Given a strongly connected  graph $G$ with $|E|$ edges and edge adjacency matrix $T_{d}$, assign to each edge $e_{i}$ a color $c_{i} \in \{c_{1}, ... , c_{|E|}\}$
so that
to a cycle of length $N$ in $G$ corresponds with an ordered sequence of $N$ colors. Assign each color in this sequence to a bead in a necklace with $N$ beads in such a manner that two adjacent colors in the sequence are assigned to adjacent beads. Then, the number of equivalence classes of nonperiodic colorings  of the necklace with $N$ beads and color matrix $T_{d}$ is $\theta(N,T_{d})$, given by (1.8).

 An important property of the directed vertex adjacency matrix $A_{d}$ of a strongly connected graph is that it is an {\it irreducible matrix}.  
Its edge adjacency matrix $T_{d}(G)$ is also irreducible in view of the fact that $T_{d}$ is the vertex adjacency matrix of the line graph of $G$ and the line graph of a strongly connected graph has the same property.

Recall that a matrix with zeros and ones is the adjacency matrix of a strongly connected graph if and only if it is irreducible [12], hence, any irreducible $n \times n$ matrix $S$ with zeros and ones as entries is a color matrix that can be interpreted as the edge adjacency matrix of a strongly connected graph $G$, and so
$\theta(N, S)$ is the number of equivalence classes of nonperiodic cycles of length $N$ in $G$ and of nonperiodic  colorings of a necklace with $N$ beads with at most $|E|=n$ colors.

The above ideas yield the following interpretation of the Ih and BL zeta functions considered in this paper. Given a color matrix  and an infinite sequence of necklaces indexed by the number of beads $N$, the Ih  and the BL functions are the functions that generate the  sequences  ($\theta(1,T), \theta(2,T), ...$) and
($\theta_{d}(1,T_{d}), \theta_{d}(2,T_{d}), ...$), respectively,  of  numbers of  classes of nonperiodic restricted colorings of the necklaces.

\subsection*{Acknowledgments}

There are several mathematicians to thank for sending me papers or references: C. Storm, V. Strehl, G. Tesler, P. Moree, J. Wojciechowski. Special thanks to Prof.
Asteroide Santana (UFSC) for help with the figures, latex commands and determinants.

\end{document}